\documentclass[oneside,english]{amsart}
\usepackage[T1]{fontenc}
\usepackage[latin9]{inputenc}
\usepackage{refstyle}
\usepackage{amsbsy}
\usepackage{amstext}
\usepackage{amsthm}
\usepackage{amssymb}
\usepackage{setspace}
\usepackage[style=numeric]{biblatex}
\makeatletter

\AtBeginDocument{\providecommand\propref[1]{\ref{prop:#1}}}
\AtBeginDocument{\providecommand\defref[1]{\ref{def:#1}}}
\AtBeginDocument{\providecommand\corref[1]{\ref{cor:#1}}}
\RS@ifundefined{subsecref}
  {\newref{subsec}{name = \RSsectxt}}
  {}
\RS@ifundefined{thmref}
  {\def\RSthmtxt{theorem~}\newref{thm}{name = \RSthmtxt}}
  {}
\RS@ifundefined{lemref}
  {\def\RSlemtxt{lemma~}\newref{lem}{name = \RSlemtxt}}
  {}

\numberwithin{equation}{section}
\numberwithin{figure}{section}
\IfFileExists{candleshape.def}{%
 \input{candleshape.def}}{}
\IfFileExists{dropshape.def}{%
 \input{dropshape.def}}{}
\IfFileExists{TeXshape.def}{%
 \input{TeXshape.def}}{}
\IfFileExists{triangleshapes.def}{%
 \input{triangleshapes.def}}{}

\theoremstyle{plain}
\newtheorem{thm}{\protect\theoremname}
\theoremstyle{definition}
\newtheorem{defn}[thm]{\protect\definitionname}
\theoremstyle{remark}
\newtheorem{rem}[thm]{\protect\remarkname}
\theoremstyle{plain}
\newtheorem{prop}[thm]{\protect\propositionname}
\theoremstyle{plain}
\newtheorem{cor}[thm]{\protect\corollaryname}

\makeatother

\usepackage{babel}
\providecommand{\corollaryname}{Corollary}
\providecommand{\definitionname}{Definition}
\providecommand{\propositionname}{Proposition}
\providecommand{\remarkname}{Remark}
\providecommand{\theoremname}{Theorem}

\begin{document}
\title{Loops of Infinite Order and Toric Foliations}
\maketitle
\begin{center}
\textit{Konstantinos Efstathiou}\\
Zu Chongzhi Center for Mathematics and Computational Science, Duke
Kunshan University
\par\end{center}

\begin{center}
\textit{Bohuan Lin$^{*}$, Holger Waalkens}\\
Bernoulli Institute for Mathematics, Computer Science and Artificial
Intelligence, University of Groningen
\par\end{center}

\begin{center}
$^{*}$e-mail: B.Lin@rug.nl
\par\end{center}
\begin{abstract}
In 2005 Dullin et al. proved in \cite{dullin2005maslov} that the
non-zero vector of Maslov indices is an eigenvector with eigenvalue
$1$ of the monodromy matrices of an integrable Hamiltonian system.
We take a close look at the geometry behind this result and extend
it to a more general context. We construct a bundle morphism defined
on the lattice bundle of an (general) integrable system, which can
be seen as a generalization of the vector of Maslov indices introduced
in \cite{dullin2005maslov}. The non-triviality of this bundle morphism
implies the existence of common eigenvectors with eigenvalue $1$
of the monodromy matrices, and gives rise to a corank $1$ toric foliation
refining the original one induced by the integrable system. Furthermore,
we show that in the case where the system has $2$ degrees of freedom,
this implies the global existence of a free $S^{1}$ action.
\end{abstract}

\section{background and motivation}

An integrable Hamiltonian system contains the following ingredients:
a $2n$ dimensional symplectic manifold $(M,\omega)$ as the phase
space, an integral map $F=(f_{1},...,f_{n}):M\rightarrow\mathbb{R}^{n}$
which is regular almost everywhere, Hamiltonian vector fields $\mathcal{X}_{i}$
for $i=1,...,n$ with $f_{i}$ as the Hamitonians and the commutativity
condition: $[\mathcal{X}_{i},\mathcal{X}_{j}]=0$. The commutativity
of the vector fields $\mathcal{X}_{i}$ induces a Hamiltonian $\mathbb{R}^{n}$
action $\Phi$ on the phase space $M$ given by $\Phi^{(t_{1},...,t_{n})}(p)=\varphi_{1}^{t_{1}}\circ...\circ\varphi_{n}^{t_{n}}(p)$
with $\varphi_{i}$ being the flow of $\mathcal{X}_{i}$ for each
$i$. The nondegeneracy of the symplectic form $\omega$ implies that
$\mathcal{X}_{1},...,\mathcal{X}_{n}$ are linearly independent wherever
$DF=(df_{1},...,df_{n})$ has full rank.

When a regular orbit of $\Phi$, i.e. an orbit consisting of regular
points of $F$, is compact, it is an $n$ dimensional Lagrangian torus.
The region $M_{0}$ consisting of compact regular orbits of $\Phi$
thereby admits a Lagrangian toric foliation $\mathfrak{F}$ with these
$n-$tori as the leaves. The Arnold-Liouville theorem shows that near
each of the leaves the foliation is nice in the sense that local action-angle
coordinates exist. However, topological and geometric obstructions
prevent in general the global existence of action-angle coordinates.
This has been investigated by Duistermaat in \cite{duistermaat1980global}.
The non-triviality of the torus bundle $F^{-1}(\mathcal{C})$ over
some loop $\mathcal{C}$ lying in the set of regular values of $F$
serves as one of these obstructions, and it is characterized by the
monodromy matrix $\mathcal{M}_{\mathcal{C}}$ of the bundle. One should
note that here for simplicity, we have assumed each level set $F^{-1}(c)$
with $c\in\mathcal{C}$ to consist of a single orbit. However, in
general, $F^{-1}(c)$ can contain more than one orbits and hence $F^{-1}(\mathcal{C})\rightarrow\mathcal{C}$
is not necessarily a torus bundle. We will deal with this problem
later in a more rigorous way.

Via the Bohr-Sommerfeld quantization of actions, their global non-existence
manifests itself in quantum mechanics as the non-existence of global
quantum numbers to assign the joint spectrum of commuting operators
that are the quantum analogues of the functions $f_{i}$. Here Maslov
indices determine whether the actions are required to be integer or
half-integer multiples of Planck's constant. From a geometric point
of view, the Maslov index of a closed curve of Lagrangian planes counts
the number of intersections of such a curve with a chosen Maslov cycle
$\Sigma$. Usually $\Sigma$ is induced by some Lagrangian vector
subbundle $\mathcal{V}$ of the tangent bundle $TM$, and the curve
of Lagrangian planes is the one induced by a curve on some Lagrangian
submanifold in the natural way. One can think of $\mathcal{V}$ as
a reference and then the Maslov index describes how the Lagrangian
planes rotate with respect to the reference along the curve. In particular,
when $M$ is the cotangent bundle of another manifold $N$, $\mathcal{V}$
is usually taken as the vertical distribution of $T^{*}N=M$.

For a regular value $c_{0}$ of $F$, a set of closed curves $\{\lambda_{1},...,\lambda_{n}\}$
on the Lagragian torus $F^{-1}(c_{0})$ can be taken such that $[\lambda_{1}],...,[\lambda_{n}]$
constitute a basis of $H_{1}\big(F^{-1}(c_{0})\big)\cong\mathbb{Z}^{n}$.
Let $a_{i}$ be the Maslov index of $\lambda_{i}$. The main result
in \cite{dullin2005maslov} shows that when it is non-zero, the vector
of Maslov indices $(a_{1},...,a_{n})$ is an eigenvector with eigenvalue
$1$ of the monodromy matrices $\mathcal{M}_{\mathcal{C}}$ for all
loops $\mathcal{C}$ of regular values of $F$ based at $c_{0}$.
As is pointed out in \cite{dullin2005maslov}, one implication of
the result is that a Hamiltonian $S^{1}$ action exists on the space
$F^{-1}(\mathcal{C})$. Moreover, it also imposes some restrictions
on the forms of the monodromy matrices.

A close look at the result in \cite{dullin2005maslov} suggests that
there should be a more general underlying structure. Namely, on the
one hand, the monodromy matrices are determined solely by the foliation
$\mathfrak{F}$. On the other hand however, the Maslov indices depend
on the choice of the Maslov cycle, or say, the Lagrangian subbundle
$\mathcal{V}$. A change of the reference may result in a change of
the triviality/non-triviality of the vector of Maslov indices, but
it affects neither on the foliation nor on the monodromy matrices.
It is then natural to ask whether we could restate the result in \cite{dullin2005maslov}
in terms of the topology of $M_{0}$ or $\mathfrak{F}$. This could
then not only give a new interpretation in different terms, but also
a generalization of the result in \cite{dullin2005maslov} to a broader
context.

It turns out that we are able to work this out with a more general
class of integrable systems which is to be defined in the next section.
These systems are called $\emph{integrable\ non-Hamiltonian\ systems}$
in the literature (e.g. see \cite{bogoyavlenskij1998extended,zung2016geometry}).
Note that the phrase `non-Hamiltonian' here means `not necessarily
Hamiltonian', and the notion also includes all integrable Hamiltonian
systems. No symplectic structure will be required for the definition
of such a system, and the $M$ is not necessarily of even dimension.

\section{basic settings and layout of the paper}

Since we are concerned with the toric foliation on the domain of compact
regular orbits of $\Phi$, we for simplicity, make some assumptions
that lead to the following definition of a regular integrable non-Hamiltonian
system with compact orbits.
\begin{defn}
\label{def:regular integrable systems of (k,n)}A regular (non-Hamiltonian)
integrable system of type $(k,n)$ with compact orbits is a triple
$(M,F,\Phi)$ where $M$ is an $k+n$ dimensional smooth manifold,
the integral map $F=(f_{1},...,f_{k}):M\rightarrow\mathbb{R}^{k}$
is a submersion and $\Phi$ is an effective $\mathbb{R}^{n}$ action
with compact orbits such that $F\circ\Phi^{(t_{1},...,t_{n})}(p)=F(p)$. 
\end{defn}

\begin{rem}
Unless otherwise stated, $M$ is assumed to be connected in this article.
\end{rem}

Since the orbits of $\Phi$ are compact, they are $n-$tori and thereby
constitute a toric foliation $\mathfrak{F}$ on $M$. As mentioned
earlier, the map $F:M\rightarrow\mathbb{R}^{n}$ is in general not
a torus bundle since a level set of $F$ may contain more than one
orbits. In order to get a torus bundle, we introduce the orbit space
$\mathcal{O}_{M}$, whose elements are the orbits of $\Phi$. The
bundle $\mathcal{L}$ of period lattices will play a central role
in our discussion. This is a bundle over the orbit space $\mathcal{O}_{M}$
where the fibers are the isotropy subgroups of the action $\Phi$.
Under Definition \ref{def:regular integrable systems of (k,n)} both
these spaces have nice structures and properties, and they will be
discussed in detail in Section $3$ and Section $4$.

The definitions of monodromy maps and monodromy matrices will be formally
introduced in Section $5$ along with their basic properties. 

The main construction and results are given in Section $6$. This
comprises a bundle epimorphism $\rho=(\rho_{1},...,\rho_{l})$ from
$\mathcal{L}$ to $\mathbb{Z}^{l}$ for some integer $l\leq n$. This
is a purely topological object associated with the system $(M,F,\Phi)$,
and as we will see in Section $7$, $\rho$ can be seen as a generalization
of the vector of Maslov indices. When $l>0$, each component of $\rho$
is a common eigenvector with eigenvalue $1$ of the (transposed) monodromy
matrices. This will have implications on the toric foliation $\mathfrak{F}$,
and in the case of $n=2$, it implies the existence a free $S^{1}$
action on $M$ which commutes with $\Phi$. The precise statements
of these results form the main results of this paper: Theorem \ref{thm:main}
and Theorem \ref{thm:main_n=00003D2 free S^1 action} in Section $6$.

In Section $7$, we revisit the case of integrable Hamiltonian systems.
We show that the vector of Maslov indices being non-zero implies the
non-triviality of $\rho$, and then show how the result in \cite{dullin2005maslov}
follows from our construction and results. Moreover, with Audin's
conjecture (which has been proved), our results lead to the existence
of a free $S^{1}$ action for integrable Hamiltonian systems on $\mathbb{R}^{4}$,
which is stated as Theorem \ref{thm:main_Hamiltonian on R^4_S^1 action}. 

The appendix is devoted to a technical proof of the local triviality
of the lattice bundle $\mathcal{L}$ and an illustration on how the
monodromy matrix $\mathcal{M}_{\mathcal{C}}$ determines the torus
bundle $F^{-1}(\mathcal{C})\rightarrow\mathcal{C}$.

\section{the orbit space $\mathcal{O}_{M}$ and the lattice bundle $\mathcal{L}$}

In this section we study the orbit space $\mathcal{O}_{M}$ and the
lattice bundle $\mathcal{L}$ associated to a regular integrable system
$(M,F,\Phi)$ of type $(k,n)$ with compact orbits.

By $\sim_{\Phi}$ we denote the equivalence relation on $M$ such
that $p\sim_{\Phi}p'$ if and only if $p$ and $p'$ are on the same
orbit of $\Phi$. The orbit space $\mathcal{O}_{M}$ is the quotient
space $M\big/\sim_{\Phi}$, that is, each element of $\mathcal{O}_{M}$
represents an orbit of $\Phi$. Then $F$ factors as $F=\bar{F}\circ q_{\Phi}$
with $q_{\Phi}:M\rightarrow\mathcal{O}_{M}=M\big/\sim_{\Phi}$ being
the quotient map and $\bar{F}:\mathcal{O}_{M}\rightarrow\mathbb{R}^{k}$.

Note that for any open set $W$ in $M$, $q_{\Phi}^{-1}\big(q_{\Phi}(W)\big)=\Phi^{\mathbb{R}^{n}}(W)$
is open. This implies that $q_{\Phi}(W)$ is open, and therefore $q_{\Phi}$
is an open map. Since each orbit is compact, $\mathcal{O}_{M}$ is
Hausdorff. Since $F$ is a submersion, for any $o\in\mathcal{O}_{M}$
and $p\in o$, there exists a local section $\sigma$ of $F$ over
some open neighbourhood $U$ of $b=F(p)$ such that $\sigma(b)=p$.
Then $V=q_{\Phi}\circ\sigma(U)=q_{\Phi}\big(\Phi^{\mathbb{R}^{n}}\circ\sigma(U)\big)$
is an open neighbourhood of $o$ and $q_{\Phi}\circ\sigma:U\rightarrow V$
is the inverse of $\bar{F}\big|_{V}:V\rightarrow U$. Therefore $\bar{F}$
is a local homeomorphism. Choose an open covering $\{V_{\alpha}\}$
of $\mathcal{O}_{M}$ such that $\bar{F}_{\alpha}=\bar{F}\big|_{V_{\alpha}}:V_{\alpha}\rightarrow U_{\alpha}$
is a homeomorphism for each $\alpha$. Then $\{(V_{\alpha},\bar{F}_{\alpha})\}$
is a smooth structure on $\mathcal{O}_{M}$ and naturally it makes
$\bar{F}$ a local diffeomorphism. Since $F=\bar{F}\circ q_{\Phi}$
and $F$ is submersive, such a smooth structure makes $q_{\Phi}$
a submersion. The discussion here yields the following proposition.
\begin{prop}
\label{prop:smooth structure of the obit space}There exists a unique
smooth structure on $\mathcal{O}_{M}$ such that $q_{\Phi}$ is a
submersion and $\bar{F}$ is a local diffeomorphism.
\end{prop}

As in the case of Hamitonian integrable systems, for any $p\in M$
the isotropy group $\mathcal{T}_{p}\subset\mathbb{R}^{n}$ of $\Phi$
at $p$ is a free Abelian group of rank $n$ \cite{cushman1997globalaspect},
and $\mathcal{T}_{p}=\mathcal{T}_{p'}$ if $p$ and $p'$ are on the
same orbit. Hence for any orbit $o\in\mathcal{O}_{M}$, define its
period lattice to be $\mathcal{L}_{o}=\mathcal{T}_{p}$ with any $p\in o$.
The lattice bundle is then defined as $\mathcal{L}=\bigsqcup_{o\in\mathcal{O}_{M}}\mathcal{L}_{o}$
and it then naturally holds $\mathcal{L}\subset\mathcal{O}_{M}\times\mathbb{R}^{n}$.
Denote by $\pi_{\mathcal{L}}$ the natural projection $\mathcal{L}\rightarrow\mathcal{O}_{M}:\ (o,T)\rightarrow o$
and endow $\mathcal{L}$ with the subspace topology inherited from
$\mathcal{O}_{M}\times\mathbb{R}^{n}$. The following basic fact holds
and its proof is inherent in the proof of the (Hamiltonian) Liouville
theorem in \cite{cushman1997globalaspect}:
\begin{thm}
\label{thm:local tivialization of the lattice bundle} For any $o\in\mathcal{O}_{M}$,
there exists a neighbourhood $U$ and smooth sections $\tau_{i}:U\rightarrow U\times\mathbb{R}^{n}$
with $i=1,...,n$, such that for each $o'\in U$, $\tau_{1}(o')$,...,$\tau_{n}(o')$
constitute a basis of $\mathcal{L}_{o'}$. As a consequence there
exists an isomorphism $\rho:\pi_{\mathcal{L}}^{-1}(U)\rightarrow U\times\mathbb{Z}^{n}$. 
\end{thm}

\begin{proof}
See Step 3 of Proof (2.1) in Chapter IX of \cite{cushman1997globalaspect}.
Note that the definition of the period lattice bundle in \cite{cushman1997globalaspect}
is different from that in this paper. In the appendix we give an argument
following the same idea but different in some technical details.
\end{proof}
We call a neighbourhood $U$ such as in the theorem above a $\mathcal{O}_{M}$-neighbourhood,
and the sections $\tau_{1}$,...,$\tau_{n}$ a smooth local basis
of $\mathcal{L}$ over $U$. Theorem \ref{thm:local tivialization of the lattice bundle}
above actually shows that $\mathcal{L}$ is a smooth submanifold of
$\mathcal{O}_{M}\times\mathbb{R}^{n}$. As a corollary of Theorem
\ref{thm:local tivialization of the lattice bundle}, we have:
\begin{cor}
The map $\pi_{\mathcal{L}}:\mathcal{L}\rightarrow\mathcal{O}_{M}$
is a covering map of $\mathcal{O}_{M}$.
\end{cor}

Since for each $o\in\mathcal{O}_{M}$, the group $\mathcal{L}_{o}$
is a subgroup of the fiber $\{o\}\times\mathbb{R}^{n}$, we can define
an equivalence relation $\sim_{\mathcal{L}}$ on $\mathcal{O}_{M}\times\mathbb{R}^{n}$
given by $(o,T)\sim_{\mathcal{L}}(o',T')$ if and only if $o=o'$
and $(o,T-T')\in\mathcal{L}_{o}$. We denote the quotient map of $\sim_{\mathcal{L}}$
by $q_{\mathcal{L}}$.

The following non-Hamiltonian version of the Liouville theorem is
a consequence of Theorem \ref{thm:local tivialization of the lattice bundle}. 
\begin{thm}
\label{thm:non-Hamiltonian Liouville theorem}\cite{zung2016geometry}
For any $o\in\mathcal{O}_{M}$, there exists some neighbourhood $U$
of $o$ in $\mathcal{O}_{M}$ with a diffeomorphism $\phi:\big(U\times\mathbb{R}^{n}\big)\big/\sim_{\mathcal{L}}\rightarrow q_{\Phi}^{-1}(U)$.
Moreover, the $\mathbb{T}^{n}$ action on $q_{\Phi}^{-1}(U)$ induced
by the natural $\mathbb{T}^{n}$ action on $\big(U\times\mathbb{R}^{n}\big)\big/\sim_{\mathcal{L}}\cong U\times\mathbb{T}^{n}$
via this diffeomorphism commutes with the $\mathbb{R}^{n}$ action
$\Phi$.
\end{thm}

\begin{proof}
Choose a local section $\sigma$ of $q_{\Phi}:M\rightarrow\mathcal{O}_{M}$
over some $\mathcal{O}_{M}$-neighbourhood $U$ of $o$. Define $\tilde{\phi}:U\times\mathbb{R}^{n}\rightarrow q_{\Phi}^{-1}(U)$
as $\tilde{\phi}(c,T)=\Phi^{T}\circ\sigma(c)$. Then $\tilde{\phi}$
factors as $\phi\circ q_{\mathcal{L}}$. It is can be checked that
$\phi$ is a homeomorphism from $\big(U\times\mathbb{R}^{n}\big)\big/\sim_{\mathcal{L}}$
to $q_{\Phi}^{-1}(U)$. 

To show that $\big(U\times\mathbb{R}^{n}\big)\big/\sim_{\mathcal{L}}$
is homeomorphic to $U\times\mathbb{T}^{n}$ one can resort to a smooth
local basis $\tau_{1}$,...,$\tau_{n}$ of $\mathcal{L}$ over $U$
and a rectifying map $rec:U\times\mathbb{R}^{n}\rightarrow U\times\mathbb{R}^{n}$
: $(c;t_{1},...,t_{n})\rightarrow(c;t_{1}\tau_{1}(c),...,t_{n}\tau_{n}(c))$.
This descends to an isomorphism $\overline{rec}$ from $U\times\mathbb{T}^{n}$
to $\big(U\times\mathbb{R}^{n}\big)\big/\sim_{\mathcal{L}}$. It can
be checked that $\overline{rec}$ is a bijective submersion, and thereby
is a diffeomorphism.
\end{proof}
\begin{cor}
\label{cor:if q_phi has a global section}If the bundle $q_{\Phi}:M\rightarrow\mathcal{O}_{M}$
admits a global section, then it is isomorphic to $\bar{\pi}_{\mathcal{O}}:\big(\mathcal{O}_{M}\times\mathbb{R}^{n}\big)\big/\sim_{\mathcal{L}}\rightarrow\mathcal{O}_{M}$.
\end{cor}

\begin{proof}
The proof of Theorem \ref{thm:non-Hamiltonian Liouville theorem}
actually guarantees that $\bar{\pi}_{\mathcal{O}}:\big(\mathcal{O}_{M}\times\mathbb{R}^{n}\big)\big/\sim_{\mathcal{L}}\rightarrow\mathcal{O}_{M}$
with $\bar{\pi}_{\mathcal{O}}([o,T])=o$ is a locally trivial torus
fibration. If $q_{\Phi}:M\rightarrow\mathcal{O}_{M}$ admits a global
section $\sigma:\mathcal{O}_{M}\rightarrow M$ then $\tilde{\phi}_{M}:\mathcal{O}_{M}\times\mathbb{R}^{n}\rightarrow M$
with $\tilde{\phi}_{M}(c,T)=\Phi^{T}\circ\sigma(c)$ for any $c\in\mathcal{O}_{M}$
factors as $\tilde{\phi}_{M}=\phi_{M}\circ q_{\mathcal{L}}$ with
$\phi_{M}$ being a bundle isomorphism.
\end{proof}

\section{the lattice bundle $\mathcal{L}$ and the sheaf $\mathcal{R}$ of
compatible $S^{1}$ actions}

The sheaf $\mathcal{R}$ of compatible $S^{1}$ actions over $\mathcal{O}_{M}$
was first introduced in \cite{zung2003symplectic} in the context
of Hamiltonian integrable systems and then in \cite{zung2004torus}
for the non-Hamiltonian case. We explain the relation between the
sheaf $\mathcal{R}$ and the lattice bundle $\mathcal{L}$ in this
section.
\begin{defn}
For any open set $U$ of $\mathcal{O}_{M}$, an $S^{1}$ action $\Theta:S^{1}\times q_{\Phi}^{-1}(U)\rightarrow q_{\Phi}^{-1}(U)$
acting on $q_{\Phi}^{-1}(U)$ is called a\textbf{ compatible $S^{1}$
action} over $U$ if it commutes with the $\mathbb{R}^{n}$ action
$\Phi$ and preserves $q_{\Phi}$ in the sense that $q_{\Phi}\circ\Theta(z,p)=q_{\Phi}(p)$
for all $z\in S^{1}$ and $p\in q_{\Phi}^{-1}(U)$. 
\end{defn}

\begin{defn}
\cite{zung2004torus} For any open set $U$ of $\mathcal{O}_{M}$,
denote by $\mathcal{R}_{U}$ the set of all compatible $S^{1}$ actions
over $U$. For any pair of open sets $U$, $V$ with $U\subset V$,
denote by $\rho_{U}^{V}$ the restriction map from $\mathcal{R}_{V}$
to $\mathcal{R}_{U}$. Then $\mathcal{R}=\big(\{\mathcal{R}_{U}\},\{\rho_{U}^{V}\big|U\subset V\}\big)$
is a sheaf of Abelian groups and it called $the\ sheaf\ of\ compatible\ S^{1}\ actions$.
\end{defn}

The following proposition characterizes the compatible $S^{1}$ actions.
\begin{prop}
\label{prop:compatible S^1 actions and sections}Suppose $\Theta:S^{1}\times q_{\Phi}^{-1}(U)\rightarrow q_{\Phi}^{-1}(U)$
is a compatible $S^{1}$ action. Then there exists a unique continuous
section $\sigma:U\rightarrow\mathcal{L}$ such that $\Theta(e^{i\cdot2\pi t},p)=\Phi^{t\cdot\sigma\circ q_{\Phi}(p)}(p)$.
Conversely, if $\sigma$ is a continuous section of $\mathcal{L}$
over $U$, then $\Theta(e^{i\cdot2\pi t},p)=\Phi^{t\cdot\sigma\circ q_{\Phi}(p)}(p)$
defines a compatible $S^{1}$ action.
\end{prop}

\begin{proof}
Let $X_{\Theta}$ be the infinitesimal generator of the flow $\psi^{t}(p)=\Theta(e^{i\cdot2\pi t},p)$
on $q_{\Phi}^{-1}(U)$, i.e. $X_{\Theta}(p)=\frac{d}{dt}\big|_{t=0}\Theta(e^{i\cdot2\pi t},p)$.
Let $\frac{\partial}{\partial t_{i}}\big|_{p}=\frac{d}{dt_{i}}\Phi^{(0,..t_{i}..,0)}(p)$.
Since $\Theta$ preserves $q_{\Phi}$, $X_{\Theta}(p)\in\ker q_{\Phi}{}_{*}=span\{\frac{\partial}{\partial t_{i}}\big|_{p},i=1,...,n\}$.
Hence there exist smooth functions $a_{1},...,a_{n}$ on $q_{\Phi}^{-1}(U)$
such that $X_{\Theta}=\Sigma_{i=1}^{n}a_{i}\frac{\partial}{\partial t_{i}}$.
Moreover, since $\Theta$ commutes with $\Phi$, the functions $a_{i}$
are invariant on each orbit, and hence they can be seen as functions
on $U$. Note that the flow $\varphi$ on $q_{\Phi}^{-1}(U)$ with
$\varphi^{t}(p)=\Phi^{t\cdot a_{1}\circ q_{\Phi}(p),...,t\cdot a_{n}\circ q_{\Phi}(p)}(p)$
has the same infinitesimal generator as $\psi$, and hence $\varphi=\psi$.
Since $\Phi^{a_{1}\circ q_{\Phi}(p),...,a_{n}\circ q_{\Phi}(p)}(p)=\Theta(e^{i\cdot2\pi},p)=p$
for all $p\in q_{\Phi}^{-1}(U)$, $\big(q_{\Phi}(p);a_{1}\circ q_{\Phi}(p),...,a_{n}\circ q_{\Phi}(p)\big)\in\mathcal{L}_{q_{\Phi}(p)}$.
Hence $\sigma:o\rightarrow\big(o;a_{1}(o),...,a_{n}(o)\big)$ is the
smooth section we are looking for. If there is another section $\sigma'$
on $U$ satisfying the relation $\Theta(e^{i\cdot2\pi t},p)=\Phi^{t\cdot\sigma'\circ q_{\Phi}(p)}(p)$
with $\sigma':o\rightarrow\big(o,a'_{1}(o),...,a'_{n}(o)\big)$, then
it holds that $X_{\Theta}=\Sigma_{i=1}^{n}a'_{i}\frac{\partial}{\partial t_{i}}$.
As a consequence, we have $a_{i}(o)=a'_{i}(o)$ for all $i=1,...,n$
and $o\in U$, and thereby $\sigma=\sigma'$.

The argument for the converse direction is straightforward.
\end{proof}
From Proposition \ref{prop:compatible S^1 actions and sections} it
follows that for any $\Theta,\Theta'\in\mathcal{R}_{U}$, the identity
$\Theta^{z}\circ\Theta'^{z'}(p)=\Theta'^{z'}\circ\Theta{}^{z}(p)$
holds, and thus there is a natural Abelian structure on $\mathcal{R}_{U}$
with the addition $\Theta+\Theta'$ given by $\big(\Theta+\Theta'\big)^{z}(p)=\Theta{}^{z}\circ\Theta'^{z}(p)$.
Accordingly, compatible $S^{1}$ actions over $U$ are in one-one
correspondence with the continuous sections of $\mathcal{L}$ over
$U$, and this implies that $\mathcal{L}$ is the associated sheaf/Etale
space of $\mathcal{R}$. The discussion here amounts to the following
corollary.
\begin{cor}
The Etale space of the sheaf $\mathcal{R}$ of compatible $S^{1}$
actions is isomorphic to the lattice bundle $\mathcal{L}$.
\end{cor}

If $\sigma:\mathcal{O}_{M}\rightarrow\mathcal{L}$ is a global continuous
section of $\mathcal{L}$, then $\sigma$ corresponds to a compatible
$S^{1}$ action $\Theta$ on $M$ by $\Theta(e^{i\cdot2\pi t},p)=\Phi^{t\cdot\sigma\circ q_{\Phi}(p)}(p)$.
When it is non-zero, this is a non-trivial $S^{1}$ action. Actually
in this case the section $\sigma$ is non-zero everywhere and the
corresponding $S^{1}$ action $\Theta$ is thereby effective. To see
this, first note that $\mathcal{O}_{M}\times\{0\}$ is closed in $\mathcal{O}_{M}\times\mathbb{R}^{n}$
and hence is closed in $\mathcal{L}$. Moreover, the fact that $\mathcal{L}$
is locally isomorphic to $U\times\mathbb{Z}^{n}$ implies that $\mathcal{O}_{M}\times\{0\}$
is also open in $\mathcal{L}$. Since we assume $M$ to be connected,
$\mathcal{O}_{M}\times\{0\}$ is exactly one connected component of
$\mathcal{L}$. Hence for any continuous section $\sigma:\mathcal{O}_{M}\rightarrow\mathcal{L}$,
as long as $\sigma(c_{0})\in\mathcal{O}_{M}\times\{0\}$, we have
$\sigma(\mathcal{O}_{M})=\mathcal{O}_{M}\times\{0\}$. As a consequence,
if $\sigma$ is non-zero somewhere, it is non-zero everywherem, and
the induced $S^{1}$ action $\Theta$ is effective.

Note that $\sigma$ being non-zero does not imply that $\Theta$ is
free. However, the existence of a non-zero section does imply the
existence of a free $S^{1}$ action compatible with $\Phi$. 
\begin{prop}
\label{prop:nonzero sections and free S^1 actions}The existence of
a global non-zero continuous section of $\mathcal{L}$ implies the
existence of a free $S^{1}$ action on $M$.
\end{prop}

\begin{proof}
First we show that for any non-zero integer $n$, the map $*n:w\mapsto n\cdot w$
defined on $\mathcal{O}_{M}\times\mathbb{R}^{n}$ induces an open
and closed map on $\mathcal{L}$. Note that $*n$ is a diffeomorphism
on $\mathcal{O}_{M}\times\mathbb{R}^{n}$, and $\mathcal{L}$ is closed
in $\mathcal{O}_{M}\times\mathbb{R}^{n}$ with $*n(\mathcal{L})\subset\mathcal{L}$.
As a consequence, $*n$ induces a closed map on $\mathcal{L}$, which
is also denoted by $*n$. For the openess, first note that for any
continuous section $\sigma$ of $\mathcal{L}$ over some open set
$U$ in $\mathcal{O}_{M}$, $\sigma(U)$ is open in $\mathcal{L}$.
Then for any $w\in\mathcal{L}$ and its neighbourhood $W$, there
exists some continuous section $\sigma$ over $U\ni c$ with $c=q_{\Phi}(w)$
such that $\sigma(U)\subset W$ and $\sigma(c)=w$. Then $*n\circ\sigma$
is also a contiuous section of $\mathcal{L}$ over the open set $U$,
and $*n\circ\sigma(U)\subset*n(W)$. Therefore $*n(W)$ is a neighbourhood
of $*n(w)$ and this implies $*n$ is an open map on $\mathcal{L}$.

Now suppose that $\sigma$ is a non-zero section. Then $\sigma(\mathcal{O}_{M})$
is a connected component of $\mathcal{L}$. Fix a point $c_{0}$ in
$\mathcal{O}_{M}$. Choose $w_{0}\in\mathcal{L}_{c_{0}}$ such that
$*n(w_{0})=\sigma(c_{0})$ for some $n\in\mathbb{Z}$, and $\mathbb{R}\cdot w_{0}\cap\mathcal{L}=\mathbb{Z}\cdot w_{0}$.
Suppose that $S_{w_{0}}\subset\mathcal{L}$ is the connected component
of $\mathcal{L}$ containing $w_{0}$. Then $*n(S_{w_{0}})$ is another
component of $\mathcal{L}$ since $*n$ is both an open and a closed
map on $\mathcal{L}$. Since $\sigma(c_{0})=*n(w_{0})\in*n(S_{w_{0}})$,
it yields $\sigma(\mathcal{O}_{M})=*n(S_{w_{0}})$, and therefore
$\sigma_{1}=\frac{1}{n}\sigma$ is also a continuous section of $\mathcal{L}$
with $\sigma_{1}(\mathcal{O}_{M})=S_{w_{0}}$. It remains to show
that for any $c\in\mathcal{O}_{M}$, the identity $\mathbb{R}\cdot\sigma_{1}(c)\cap\mathcal{L}=\mathbb{Z}\cdot\sigma_{1}(c)$
holds. That is, $\sigma_{1}(c)$ is a generator of the subgroup $\mathbb{R}\cdot\sigma_{1}(c)\cap\mathcal{L}$
for each $c\in$$\mathcal{O}_{M}$.

Suppose that $w_{1}\in\mathcal{L}_{c}$ is a generator of $\mathbb{R}\cdot\sigma_{1}(c)\cap\mathcal{L}$
. Then there exists some integer $n_{1}$ such that $n_{1}w_{1}=\sigma_{1}(c)$.
Repeating the argument above yields another section $\sigma_{2}=\frac{1}{n_{1}}\sigma_{1}$.
Since $\sigma_{1}(c_{0})=w_{0}$ and $\mathbb{R}\cdot w_{0}\cap\mathcal{L}=\mathbb{Z}\cdot w_{0}$,
$\sigma_{2}(c_{0})=\frac{1}{n_{1}}\sigma_{1}(c_{0})=\frac{1}{n_{1}}w_{0}\in\mathbb{R}\cdot w_{0}\cap\mathcal{L}=\mathbb{Z}\cdot w_{0}$.
This implies $n_{1}=\pm1$. Hence for any $c\in\mathcal{O}_{M}$,
$\mathbb{R}\cdot\sigma_{1}(c)\cap\mathcal{L}=\mathbb{Z}\cdot\sigma_{1}(c)$,
and therefore the $S^{1}$ action $\Theta$ on $M$ defined by $\Theta(e^{i\cdot2\pi t},p)=\Phi^{t\cdot\sigma_{1}\circ q_{\Phi}(p)}(p)$
is free.
\end{proof}

\section{monodromy maps}

Consider a loop $\mathcal{C}$ in the orbit space $\mathcal{O}_{M}$
with a fixed point $c_{0}\in\mathcal{C}$. For convenience, we view
$\mathcal{C}$ both as a subset $\mathcal{C}\subset\mathcal{O}_{M}$
and a fixed parametrization $\mathcal{C}:[0,1]\rightarrow\mathcal{O}_{M}$
with $\mathcal{C}(0)=\mathcal{C}(1)=c_{0}$. Let $M_{\mathcal{C}}=q_{\Phi}^{-1}(\mathcal{C})$
and let $\mathcal{L}_{\mathcal{C}}=\pi_{\mathcal{L}}^{-1}(\mathcal{C})$.
Then $q_{\mathcal{C}}=q_{\Phi}\big|_{\mathcal{C}}:M_{\mathcal{C}}\rightarrow\mathcal{C}$
is a locally trivial torus fibration over $\mathcal{C}$. Since $q_{\mathcal{C}}$
is a locally trivial fibration over a loop with connected fibers,
it always admits a section over $\mathcal{C}$. 

According to Corollary $\ref{cor:if q_phi has a global section}$,
$q_{\mathcal{C}}:M_{\mathcal{C}}\rightarrow\mathcal{C}$ is isomorphic
to $\bar{\pi}_{\mathcal{C}}:\big(\mathcal{C}\times\mathbb{R}^{n}\big)\big/\sim_{\mathcal{L}}\rightarrow\mathcal{C}$.
Recall that $\mathcal{L}\rightarrow\mathcal{O}_{M}$ is a covering,
and denote by $\mathcal{M}_{\mathcal{C}}$ the monodromy action of
$\mathcal{C}$ on $\mathcal{L}_{c_{0}}$. $\mathcal{M}_{\mathcal{C}}$
is actually an isomorphism on the lattice $\mathcal{L}_{c_{0}}$.
To see this, suppose $\gamma_{1}$ and $\gamma_{2}$ are the lifts
of $\mathcal{C}$ with base points $v_{1}$ and $v_{2}$ , with, $v_{1},v_{2}\in\mathcal{L}_{c_{0}}$.
$\gamma_{1}+\gamma_{2}$ is the lift of $\mathcal{C}$ at $v_{1}+v_{2}$,
which implies that $\mathcal{M}_{\mathcal{C}}$ is a group homomorphism.
It can be checked that $\mathcal{M}_{\mathcal{C}}$ is bijective.
By fixing a basis $\bar{w}=(u_{1},...,u_{n})$ of $\mathcal{L}_{c_{0}}$,
$\mathcal{M}_{\mathcal{C}}$ is represented by some element $\mathcal{M}_{\mathcal{C},\bar{w}}\in SL(n,\mathbb{Z})$.
Note that $\mathcal{M}_{\mathcal{C},\bar{w}}\in GL(n,\mathbb{Z})$.
To see that $\det\mathcal{M}_{\mathcal{C},\bar{w}}=1$, we only need
to show that $\det\mathcal{M}_{\mathcal{C},\bar{w}}>0$. For each
$i\in\{1,...,n\}$, let $\tau_{i}:[0,1]\rightarrow\mathcal{L}\subset\mathcal{O}_{M}\times\mathbb{R}^{n}$
be the lift of $\mathcal{C}$ with $\tau_{i}(0)=u_{i}$. It holds
that $\tau_{i}(s)=\big(\mathcal{C}_{t},\tau'_{i}(s)\big)$ with $\tau'_{i}:[0,1]\rightarrow\mathbb{R}^{n}$
being continuous for each $i\in\{1,...,n\}$. Consequently $\det[\tau_{1}'(s),...,\tau_{n}'(s)]$
is continuous and non-zero everywhere with respect to $s$ and hence
its sign does not change. By the definition of $\mathcal{M}_{\mathcal{C},\bar{w}}$,
it holds that $[\tau_{1}'(1),...,\tau_{n}'(1)]=[\tau_{1}'(0),...,\tau_{n}'(0)]\cdot\mathcal{M}_{\mathcal{C},\bar{w}}$
and thereby $\det\mathcal{M}_{\mathcal{C},\bar{w}}>0$.

We call $\mathcal{M}_{\mathcal{C}}$ the monodromy map associated
to $\mathcal{C}$ (or, of the fibration $q_{\mathcal{C}}:M_{\mathcal{C}}\rightarrow\mathcal{C}$),
and $\mathcal{M}_{\mathcal{C},\bar{w}}$ the monodromy matrix with
respect to $w$. The monodromy map $\mathcal{M}_{\mathcal{C}}$ determines
the structure of $q_{\mathcal{C}}:M_{\mathcal{C}}\rightarrow\mathcal{C}$
(see Appendix 8.2). 

\section{the main construction}

In this section we construct a bundle morphism $\rho:\mathcal{L}\rightarrow\mathbb{Z}^{l}$
with $l$ being an integer and show the main results of this article. 

Let $o$ be an element in $\mathcal{O}_{M}$. For $T\in\mathcal{L}_{o}$
and $x\in q_{\Phi}^{-1}(o)$, let $\lambda_{T,x}:[0,1]\rightarrow q_{\Phi}^{-1}(o)$
be the closed path on the torus $q_{\Phi}^{-1}(o)$ with $\lambda_{T,x}(s)=\Phi^{s\cdot T}(x)$.
Denote by $Tor_{H_{1}(M)}$ the torsion subgroup of $H_{1}(M)$. Define
$\rho_{o}:\mathcal{L}_{o}\rightarrow H_{1}(M)\big/Tor_{H_{1}(M)}$
by assigning to each $T\in\mathcal{L}_{o}$ the element $[\lambda_{T,x}]$
in $H_{1}(M)\big/Tor_{H_{1}(M)}$. Note that such an assignment is
independent of $x$ since for any $x,x'$ on $q_{\Phi}^{-1}(o)$,
$\lambda_{T,x}$ and $\lambda_{T,x'}$ are homotopic on $q_{\Phi}^{-1}(o)$.
The map $\rho_{o}$ is a homomorphism between Abelian groups. 
\begin{defn}
The bundle morphism $\rho:\mathcal{L}\rightarrow H_{1}(M)\big/Tor_{H_{1}(M)}$
is the one satisfying the identity $\rho\big|_{\mathcal{L}_{o}}=\rho_{o}$
for each $o\in\mathcal{O}_{M}$.
\end{defn}

We show that the value of $\rho$ is invariant under parallel translation
on $\mathcal{L}$, and as a consequence, $\rho$ can be seen as a
continuous bundle epimorphism from $\mathcal{L}$ to $\mathbb{Z}^{l}$
with $l\leq n$.
\begin{thm}
For any path $\tilde{\gamma}:[0,1]\rightarrow\mathcal{L}$, $\rho\circ\tilde{\gamma}(t)=\rho\circ\tilde{\gamma}(0)$
for all $t\in[0,1]$. As a consequence, $\rho$ is locally constant,
and $\text{Im}\rho=\text{Im}\rho_{o}$ is a free Abelian group with
rank $l\leq n$.
\end{thm}

\begin{proof}
Let $\gamma$ be the path in $\mathcal{O}_{M}$ defined as $\gamma=\pi_{\mathcal{L}}\circ\tilde{\gamma}$.
Take a section $\sigma$ of $q_{\Phi}:M\rightarrow\mathcal{O}_{M}$
over $\gamma$. Then $(t,\cdot)\mapsto\lambda_{\tilde{\gamma}(t),\sigma\circ\gamma(t)}(\cdot)$
gives a homotopy in $M$ and hence $\rho\circ\tilde{\gamma}(t)=[\lambda_{\tilde{\gamma}(t),\sigma\circ\gamma(t)}]$
remains invariant in $H_{1}(M)\big/Tor_{H_{1}(M)}$. Since $\mathcal{O}_{M}$
is path-connected, for any $o,o'\in\mathcal{O}_{M}$, $\mathcal{L}_{o'}$
can be obtained via the parallel translation of $\mathcal{L}_{o}$
along some path connecting $o$ and $o'$. Therefore $\text{Im}\rho_{o}=\text{Im}\rho_{o'}$
and it is a finitely generated subgroup in $H_{1}(M)\big/Tor_{H_{1}(M)}$
and thereby a free Abelian group with rank no larger than that of
$\mathcal{L}_{o}$. Finally, $\rho$ is locally constant since $\mathcal{L}$
is locally path-connected.
\end{proof}
\begin{cor}
The map $\rho$ is a bundle epimorphism from $\mathcal{L}$ to $\mathbb{Z}^{l}$
with $l\leq n$. For any path $\gamma$ in $\mathcal{O}_{M}$ with
$\mathcal{M}_{\gamma}$ being the associated monodromy map, $\rho\circ\mathcal{M}_{\gamma}=\rho$
on $\mathcal{L}_{\gamma(0)}$. 
\end{cor}

Suppose $l>0$, i.e. $\rho$ is non-trivial. Then $\rho=(\rho_{1},...,\rho_{l})$
with $\rho_{i}:\mathcal{L}\rightarrow\mathbb{Z}$ being locally constant
and $\rho_{i}\circ\mathcal{M}_{\gamma}=\rho_{i}$. This means the
linear functionals $\rho_{i}\big|_{\mathcal{L}_{\gamma(0)}}$ are
eigenvectors with eigenvalue $1$ of the transpose of $\mathcal{M}_{\gamma}$.
This gives a descending chain of lattice subbundles $\mathcal{L}\supset\ker\rho_{1}\supset\ker(\rho_{1},\rho_{2})\supset...\supset\ker(\rho_{1},...,\rho_{l})$.
Since $\mathbb{Z}$ is free, on each fiber $\mathcal{L}_{o}$ the
following short exact sequence splits:
\[
0\rightarrow\ker\rho_{1}\big|_{\mathcal{L}_{o}}\rightarrow\mathcal{L}_{o}\xrightarrow{\rho_{1}}\mathbb{Z}\rightarrow0,
\]
and then $\ker\rho_{1}\big|_{\mathcal{L}_{o}}$ is isomorphic to $\mathbb{Z}^{n-1}$. 
\begin{thm}
The sublattice bundle $\ker\rho_{1}$ is a smooth lattice subbundle
of $\mathcal{L}$. Moreover, $\ker\rho_{1}$ locally splits $\mathcal{L}$.
To be precise, for any point $o\in\mathcal{O}_{M}$, there exists
some neighbourhood $U$ such that $\mathcal{L}\big|_{U}=\ker\rho_{1}\big|_{U}\oplus\mathcal{L}_{U}^{''}$
with $\mathcal{L}_{U}^{''}$ being some sublattice bundle of $\mathcal{L}\big|_{U}$
over $U$.
\end{thm}

\begin{proof}
Since $\rho_{1}$ is locally constant on $\mathcal{L}$, it is contant
on each of the connected components. Hence $\ker\rho_{1}$ consists
of several connected components of $\mathcal{L}$ and is a submanifold
of $\mathcal{L}$.

For any $c\in\mathcal{O}_{M}$, we construct a local trivialization
of $\ker\rho_{1}$ in the vicinity of $c$. Let $U\ni c$ be a connected
neighbourhood over which $\mathcal{L}$ admits a local trivialization.
Then for each connected component $\mathcal{S}$ of $\mathcal{L}\big|_{U}$,
there exists a section $\sigma:U\rightarrow\mathcal{L}$ such that
$\mathcal{S}=\sigma(U)$. Due to the connectedness of $U$, we have
$\sigma(U)\subset\ker\rho_{1}$ if and only if $\sigma(U)\cap\ker\rho_{1}$
is nonempty. Fix a basis $w_{1},...,w_{n-1}$ of $\ker\rho_{1}\big|_{\mathcal{L}_{c}}$
with local sections $\sigma_{1},...,\sigma_{n-1}$ over $U$ such
that $\sigma_{i}(c)=w_{i}$. For each $c'\in U$, the linear independence
of $\sigma_{1}(c'),...,\sigma_{n-1}(c')$ follows from the fact that
the zero section $U\times\{\boldsymbol{0}\}$ is a component of $\ker\mathcal{L}\big|_{U}$. 

Now we show that for any $w'\in\ker\rho_{1}\big|_{\mathcal{L}_{c'}}$
, there exist integers $k_{1},...,k_{n-1}$ such that $\Sigma_{i}k_{i}\sigma_{i}(c')=w'$.
Let $\sigma':U\rightarrow\mathcal{L}$ be the section with $w'=\sigma'(c')$.
Then it holds that $\sigma'(U)\subset\ker\rho_{1}$ due to the connectness
of $U$, and in particular, $\sigma'(c)\in\ker\rho_{1}\big|_{\mathcal{L}_{c}}$.
Then there exist integers $k_{1},...,k_{n-1}$ such that $\sum_{i=1}^{n}k_{i}\sigma_{i}(c)=\sigma'(c)\in\mathcal{L}'_{c}$.
Since $\tau=\sum_{i}k_{i}\sigma_{i}$ is also a continuous section
of $\mathcal{L}$ over $U$, it yields $\tau(U)\subset\ker\rho_{1}$
and thereby $\tau=\sigma'$. Then $w'=\sigma'(c')=\sum_{i}k_{i}\sigma_{i}(c')$
with $k_{i}$ integers. Hence $\{\sigma_{1},...,\sigma_{n-1}\}$ gives
rise to a local trivialization of $\ker\rho_{1}$ over $U$.

To see that $\ker\rho_{1}$ locally splits $\mathcal{L}$, first note
that for every $c\in\mathcal{O}_{M}$, $\ker\rho_{1}\big|_{\mathcal{L}_{c}}$
splits $\mathcal{L}_{c}$. Then it holds that $\mathcal{L}_{c}=\mathcal{L}'_{c}\oplus\mathbb{Z}\cdot v_{c}$
for some $v_{c}\in\mathcal{L}_{c}$. Let $\sigma$ be a section over
$U$ such that $\sigma(c)=v_{c}$. For any $c'\in U$ and $w'\in\mathcal{L}{}_{c'}$,
there is a section $\sigma'$ of $\mathcal{L}$ over $U$ such that
$\sigma'(c')=w'$. Since $\sigma'(c)\in\mathcal{L}_{c}$ and $\sigma_{1}(c),...,\sigma_{n-1}(c),\sigma(c)$
constitute a basis of $\mathcal{L}_{c}$, there exist integers $b_{1},...,b_{n-1},b$
such that $\sigma'(c)=\sum_{i=1}^{n-1}b_{i}\sigma_{i}(c)+b\sigma(c)$.
It then holds $\big(\sum_{i=1}^{n}b_{i}\sigma_{i}\big)(U)=\sigma'(U)$
since both of them are connected components of $\mathcal{L\big|}_{U}$.
Hence $\sum_{i=1}^{n-1}b_{i}\sigma_{i}(c')+b\sigma(c')=\sigma'(c')=w'$.
As a result, $\mathcal{L}\big|_{U}=\ker\rho_{1}\big|_{U}\oplus\mathcal{L}_{U}^{''}$
with $\mathcal{L}^{''}{}_{U}=\mathbb{Z}\cdot\sigma(U)$.
\end{proof}
\begin{cor}
$\ker\rho_{1}$ gives rise to an $(n-1)-$ toric foliation $\mathfrak{F}^{(1)}$
refining $\mathfrak{F}$. More precisely, for any $b\in\mathcal{O}_{M}$,
there exists a neighbourhood $U$ of $b$ with a local trivialization
$\bar{\psi}:U\times S^{1}\times\mathbb{T}^{n-1}\rightarrow q_{\Phi}^{-1}(U)$
such that for any $(u,z)\in U\times S^{1}$, $\bar{\psi}\big((u,z)\times\mathbb{T}^{n-1}\big)$
is a leaf of $\mathfrak{F}^{(1)}$.
\end{cor}

\begin{proof}
Denote by $V_{o}$ the $n-1$ dimensional subspace of $\{o\}\times\mathbb{R}^{n}$
spanned by $\ker\rho_{1}\big|_{\mathcal{L}{}_{o}}$. Then for each
$p\in q_{\Phi}^{-1}(o)$, $\mathfrak{T}_{p}=\Phi^{V_{o}}(p)$ is an
$(n-1)-$ torus. To see $\mathfrak{F}^{(1)}=\{\mathfrak{T}_{p}\big|p\in M\}$
is a foliation, it only needs to be checked that there are local flat
charts everywhere on $M$ \cite{lee2013distributions}. 

For any $p\in M$ with $o=q_{\Phi}(p)$, choose a neighourhood $U$
over which it admits a section $\tilde{\sigma}$ of $q_{\Phi}$ with
$\tilde{\sigma}(o)=p$ and trivializations of $\mathcal{L}_{U}$.
Note that each component of $\mathcal{L}_{U}$ takes the form $\sigma(U)$
with $\sigma$ being a continuous section of $\mathcal{L}$ over $U$.
Choose $w_{1},...,w_{n-1}\in\ker\rho_{1}\big|_{\mathcal{L}{}_{o}}$
and $v\in\mathcal{L}_{o}$ such that $\{w_{1},...,w_{n-1},v\}$ forms
a basis of $\mathcal{L}_{o}$ and $\{w_{1},...,w_{n-1}\}$ forms a
basis of $\ker\rho_{1}\big|_{\mathcal{L}{}_{o}}$. Let $\sigma_{1},...,\sigma_{n-1},\tau$
be sections of $\mathcal{L}$ over $U$ such that $\sigma_{i}(o)=w_{i}$
and $\tau(o)=v$. Then for each $c\in U$, $\{\sigma_{1}(c),...,\sigma_{n-1}(c),\tau(c)\}$
is a basis for $\mathcal{L}_{c}$ with $\{\sigma_{1}(c),...,\sigma_{n-1}(c)\}$
being a basis for $\ker\rho_{1}\big|_{\mathcal{L}_{c}}$.

Define a map $\psi$: $U\times\mathbb{R}^{n}\rightarrow M$ by
\[
(u;s,t_{1},..,t_{n-1})\rightarrow\Phi^{s\cdot\tau(c)+t_{1}\cdot\sigma_{1}(c)+...+t_{n-1}\cdot\sigma_{n-1}(c)}\circ\tilde{\sigma}(c).
\]
It factors as $\psi=\bar{\psi}\circ q_{\mathbb{Z}^{n}}$ with $q_{\mathbb{Z}^{n}}:U\times\mathbb{R}^{n}\rightarrow U\times\mathbb{T}^{n}\cong U\times(\mathbb{R}^{n}\big/\mathbb{Z}^{n})$
the quotient map and $\bar{\psi}=U\times\mathbb{T}^{n}\rightarrow W=\psi(U\times\mathbb{R}^{n})$
a diffeomorphism. For any $p'\in\psi(u',s',t'_{1},...,t'_{n-1})$,
\[
\bar{\psi}^{-1}(\mathfrak{T}_{p'}\cap W)=\bar{\psi}^{-1}(\mathfrak{T}_{p'})=\{(u,\bar{s},\bar{t}_{1},...,\bar{t}_{n-1})\in U\times\mathbb{T}^{n}\big|u=u',\bar{s}=e^{i2\pi\cdot s'}\}.
\]
\end{proof}
Note that $\rho$ is an epimorphism and therefore $\ker(\rho_{1},\rho_{2})$
is a lattice subbundle of $\ker\rho_{1}$ and has corank $1$ in $\ker\rho_{1}$.
Consequently it induces an $(n-2)-$ toric foliation $\mathfrak{F}^{(2)}$
on $M$ that refines $\mathfrak{F}^{(1)}$. This process can be iterated
for all $i=1,...,l$ and then the following theorem holds.
\begin{thm}
\label{thm:main}Let $(M,F,\Phi)$ be a regular integrable system
of type $(k,n)$ with compact orbits. Suppose that $M$ (or equivalently,
its orbit space $\mathcal{O}_{M}$) is connected. If for some point
$c_{0}\in\mathcal{O}_{M}$, there exists a loop $\lambda$ on $q_{\Phi}^{-1}(c_{0})$
such that $[\lambda]$ has infinite order in $H_{1}(M)$, then the
lattice bundle $\mathcal{L}$ of the system has a series of sublattice
bundles $\mathcal{L}^{(1)}\supset...\supset\mathcal{L}^{(l)}$ for
some positive integer $l$ no more than $n$ such that, $\mathcal{L}^{(i)}=\ker(\rho_{1},...,\rho_{l})$
has rank $n-i$ for each $i\in\{1,...,l\}$, $\mathcal{L}^{(i+1)}$
locally splits $\mathcal{L}^{(i)}$ for each $i\in\{1,...,l-1\}$
and $\mathcal{L}^{(1)}$ locally splits $\mathcal{L}$. As a consequence,
such a sequence of sublattice bundles gives rise to a sequence of
toric foliations $\mathfrak{F}^{(1)},...,\mathfrak{F}^{(l)}$ such
that $\mathfrak{F}^{(i+1)}$ refines $\mathfrak{F}^{(i)}$ for $i=1,...,l-1$
and $\mathfrak{F}^{(1)}$ refines the fibration $\mathfrak{F}$. For
each $i\in\{1,...,l\}$, the leaves of $\mathfrak{F}^{(i)}$ are $(n-i)-$
tori.
\end{thm}

In the case $n=2$, when $\rho:\mathcal{L}\rightarrow H_{1}(M)\big/Tor_{H_{1}(M)}$
is non-trivial, $\text{Im\ensuremath{\rho}}$ is isomorphic to either
$\mathbb{Z}$ or $\mathbb{Z}^{2}$, and then we at least get one non-zero
linear functional $\rho_{1}:\mathcal{L}\rightarrow\mathbb{Z}$. Then
$\ker\rho_{1}$ is a lattice subbundle with rank $1$. Actually $\ker\rho_{1}\cong\mathcal{O}_{M}\times\mathbb{Z}$
and thus by Proposition \ref{prop:nonzero sections and free S^1 actions}
there is a free $S^{1}$ action on $M$.
\begin{thm}
\label{thm:main_n=00003D2 free S^1 action}In the case $n=2$, if
$\rho:\mathcal{L}\rightarrow H_{1}(M)\big/Tor_{H_{1}(M)}$ is non-trivial,
there exists a compatible free $S^{1}$ action on $M$. 
\end{thm}

\begin{proof}
We only need to show that $\ker\rho_{1}\cong\mathcal{O}_{M}\times\mathbb{Z}$.
Fix a point $c_{0}$ in $\mathcal{O}_{M}$. For any other point $o\in\mathcal{O}_{M}$,
$\ker\rho_{1}\big|_{\mathcal{L}_{o}}$ can be obtained via the parallel
transport of $\ker\rho_{1}\big|_{\mathcal{L}_{c_{0}}}$ along any
path joining $c_{0}$ and $o$. Choose a(n) (ordered) basis $w=(u,v)$
of $\mathcal{L}_{c_{0}}$ and let $z\in\mathbb{Z}^{2}$ with $z=\big(\rho_{1}(u),\rho_{1}(v)\big)$.
For any closed path $\gamma$ in $\mathcal{O}_{M}$ with $\gamma(0)=\gamma(1)=c_{0}$,
$z\cdot\mathcal{M}_{\gamma,w}=z$ with $\mathcal{M}_{\gamma,w}$ being
the monodromy matrices. Hence both of the eigenvalues of $\mathcal{M}_{\gamma,w}$
equal to $1$. Since $\ker\rho_{1}\big|_{\mathcal{L}_{c_{0}}}$ is
a one dimensional invariant space of the monodromy map $\mathcal{M}_{\gamma}$,
it holds that $\mathcal{M}_{\gamma}\cdot T=T$ for all $T\in\ker\rho_{1}\big|_{\mathcal{L}_{c_{0}}}$.
Due to the arbitrariness of $\gamma$, this implies the triviality
of the bundle $\ker\rho_{1}\rightarrow\mathcal{O}_{M}$.
\end{proof}

\section{maslov indices, monodromy matrices and toric foliations with corank
1}

Now we restrict to the Hamiltonian context and show how the results
obtained in the previous section are related to the work in \cite{dullin2005maslov}. 

Consider an integrable Hamiltonian system $(M,\omega,F)$ with $F=(f_{1},...,f_{n})$
being the integral map. Let $\mathcal{X}_{i}$ be the Hamiltonian
vector fields with $df_{i}(\cdot)=\omega(\mathcal{X}_{i},\cdot)$.
We assume that the orbits of the Hamiltonian $\mathbb{R}^{n}$ action
$\Phi$ are all compact. Denote by $\Lambda_{M}$ the bundle of Lagrangian
Grassmanians of $M$. 

Fix a complex structure $\mathcal{J}$ compatible with the symplectic
form $\omega$, and let $g_{\mathcal{J}}(\cdot,*)=\omega(\mathcal{J}\cdot,*)$
be the compatible Riemannian structure. Recall that the Hamiltonian
vector fields $\mathcal{X}_{1},...,\mathcal{X}_{n}$ are independent
everywhere by the assumption we made and $span_{\mathbb{R}}\{\mathcal{X}_{1},...,\mathcal{X}_{n}\}$
is a Lagrangian distribution. Apply the Gram-Schmidt process to obtain
$\mathcal{X}'_{1},...,\mathcal{X}'_{n}$ that are orthonormal with
repect to $g_{\mathcal{J}}$. Note that $span_{\mathbb{R}}\{\mathcal{X}'_{1},...,\mathcal{X}'_{n}\}=span_{\mathbb{R}}\{\mathcal{X}_{1},...,\mathcal{X}_{n}\}$
and it is Lagrangian. Then,
\[
\{\mathcal{X}'_{1},...,\mathcal{X}'_{n},\mathcal{J}(\mathcal{X}'_{1}),...,\mathcal{J}(\mathcal{X}'_{n})\}
\]
 is a globally defined symplectic frame for the tangent bundle $TM$. 

Recall that for a symplectic manifold $(M,\omega)$, the symplectic
form $\omega$ is a symplectic bilinear form on its tangent bundle,
and $(TM,\omega,M)$, as a vector bundle endowed with a symplectic
bilinear form, is a symplectic vector bundle (see page 79 in \cite{mcduff2017introduction}).
Then the argument above shows that the symplectic vector bundle $(TM,\omega,M)$
is isomorphic to the trivial symplectic vector bundle $(M\times\mathbb{R}^{2n},\omega_{0},M)$
with $\omega_{0}\big|_{\{p\}\times\mathbb{R}^{2n}}=dx_{i}\wedge dy_{i}$
being the standard symplectic bilinear form on $\{p\}\times\mathbb{R}^{2n}$.
Therefore, there exists an isomorphism $\Lambda_{M}\cong M\times\mathbb{U}(n)\big/\mathbb{O}(n)$
associated to such a trivialization. Let $\Lambda_{M_{\mathcal{C}}}\rightarrow M_{\mathcal{C}}$
denote the restriction of $\Lambda_{M}$ to $M_{\mathcal{C}}$ and
then $\Lambda_{M_{\mathcal{C}}}\cong M_{\mathcal{C}}\times\mathbb{U}(n)\big/\mathbb{O}(n)$.
Let $E:M\rightarrow\Lambda_{M}$ be a section, i.e. a Lagrangian vector
bundle over $M$. Then for any $p\in M$, $E(p)=[E_{p}]\in\mathbb{U}(n)\big/\mathbb{O}(n)$
with $E_{p}\in\mathbb{U}(n)$ being a representative.

For any section $E$, define $m_{E}:\Lambda_{M}\cong M\times\mathbb{U}(n)\big/\mathbb{O}(n)\rightarrow S^{1}$
with $m_{E}(p,[A])=\big(det_{\mathbb{C}}(A\cdot E_{p}^{-1})\big)^{2}$.
Note that for any $p\in M$, $span_{\mathbb{R}}\big(\mathcal{X}_{1}(p),...,\mathcal{X}_{n}(p)\big)\subset\Lambda_{M}\big|_{p}$
and hence there is a natural embedding $\tilde{l}$ of $M$ into $\Lambda_{M}$
by $p\rightarrow span_{\mathbb{R}}\big(\mathcal{X}_{1}(p),...,\mathcal{X}_{n}(p)\big)$. 
\begin{defn}
The $Maslov\ map$ for the integrable system $(M,F,\Phi)$ is the
map $\tilde{m}_{E}=m_{E}\circ\tilde{l}$. 
\end{defn}

\begin{defn}
\label{def:The-Maslov-index} The $Maslov\ index$ $\mu_{\lambda}$
of a loop $\lambda$: $S^{1}\rightarrow M$ in the integrable system
$(M,F,\Phi)$ with respect to $E$ is the degree of the map $\tilde{m}_{E}\circ\lambda$.
\end{defn}

Definition \ref{def:The-Maslov-index} is slightly different from
the usual definition of the Maslov indices in that we do not require
the loop $\lambda$ to be on some Lagrangian submanifold, since the
integrable system already prescribes to each point in $M$ a Lagrangian
subspace. Yet this definition is consistent with the usual one when
$\lambda$ lies on a Lagrangian submanifold.

Note that for any fixed point $x_{0}\in M$, the homomorphism $\pi_{1}(M,x_{0})\ni[\lambda]\mapsto\mu_{\lambda}\in\mathbb{Z}$
factors through $H_{1}(M)\big/Tor_{H_{1}(M)}$ since $\mathbb{Z}$
is Abelian and torsionless. That is, there is a homomorphism $\mu:H_{1}(M)\big/Tor_{H_{1}(M)}\rightarrow\mathbb{Z}$
such that for any loop $\lambda$ in $M$, $\mu([\lambda])=\mu_{\lambda}$.
The composition $\mu\circ\rho:\mathcal{L}\rightarrow\mathbb{Z}$ is
then a bundle morphism. Note that for $c\in\mathcal{O}_{M}$, $\mu\circ\rho(\mathcal{L}_{c})=n_{c}\cdot\mathbb{Z}$
for some non-negative integer $n_{c}$. $n_{c}$ is called the minimal
Maslov number of the Lagrangian torus $q^{-1}(c)$. If the minimal
Maslov number on some Lagrangian torus $q_{\Phi}^{-1}(c_{0})$ is
non-zero, then $\mu\circ\rho$ is non-zero, implying $\rho$ to be
non-trivial. Then $\mu\circ\rho\big|_{\mathcal{L}_{c_{0}}}$ is a
common eigenvector with eigenvalue $1$ for the transpose of the monodromy
maps $\mathcal{M}_{\gamma}$ of all the loops $\gamma$ in $\mathcal{O}_{M}$
with $c_{0}\in\gamma$. For a basis $(u_{1},...,u_{n})$ of $\mathcal{L}_{c_{0}}$,
$\big(\mu\circ\rho(u_{1}),...,\mu\circ\rho(u_{n})\big)$ gives the
corresponding vector of Maslov indices. In this way, it yields Theorem
$1$ of \cite{dullin2005maslov}. 

We conclude this article with a result for integrable Hamiltonian
systems in $\mathbb{R}^{4}$.

Audin's conjecture \cite{audin1988fibres} asserts that for any Lagrangian
torus in $(\mathbb{R}^{2n},dx_{1}\wedge dy_{1}+\cdots+dx_{n}\wedge dy_{n})$,
the minimal Maslov number is $2$ and this has been confirmed in \cite{cieliebak2018punctured}
following a series of partial results \cite{viterbo1990new,buhovsky2010maslov,fukaya2006application,damianfloer}
(see page 118 in \cite{mcduff2017introduction} for a brief introduction
to the results obtained in these papers). Combining this result with
Theorem \ref{thm:main_n=00003D2 free S^1 action} yields:
\begin{thm}
\label{thm:main_Hamiltonian on R^4_S^1 action}Suppose that $(\mathbb{R}^{4},\omega,F)$
is a Hamiltonian integrable system (not necessarily regular) with
integral map $F:\mathbb{R}^{4}\rightarrow\mathbb{R}^{2}$. Denote
by $\mathbb{R}_{reg}^{4}$ the set of regular points of $F$. If $M_{0}$
is a connected domain of $\mathbb{R}_{reg}^{4}$ within which the
orbits of the Hamiltonian $\mathbb{R}^{2}$ action are compact, then
there exists a free Hamiltonian $S^{1}$ action on $M_{0}$.
\end{thm}

\section{Appendix }

\subsection{Local trivialization of the lattice bundle}

Let $(M,F,\Phi)$ be a regular integrable system of type $(k,n)$
with compact $\Phi-$ orbits. We give a detailed argument for Theorem
\ref{thm:local tivialization of the lattice bundle} in Section 3.
For convenience, we state Theorem \ref{thm:local tivialization of the lattice bundle}
again.

\textbf{Theorem \ref{thm:local tivialization of the lattice bundle}.}
For any $o\in\mathcal{O}_{M}$, there exists some neighbourhood $U_{o}$
over which there are smooth sections $\tau_{i}:U\rightarrow U\times\mathbb{R}^{n}$
with $i=1,...,n$, such that for each $c\in U$, $\tau_{1}(c)$,...,$\tau_{n}(c)$
constitute a basis of $\mathcal{L}_{c}$. As a consequence, there
exists an isomorphism (algebraically and topologically) $\beta:\pi_{\mathcal{L}}^{-1}(U)\rightarrow U\times\mathbb{Z}^{n}$. 
\begin{proof}
As is shown in Proposition \propref{smooth structure of the obit space},
$q_{\Phi}:M\rightarrow\mathcal{O}_{M}$ is a submersion. Hence there
exists some smooth section $\sigma:U\rightarrow M$ over some open
neighbourhood $U$ of $o$. Define $\Psi:U\times\mathbb{R}^{n}\rightarrow M$
as $\Psi(b,T)\rightarrow\Phi^{T}\circ\sigma(b)$. According to Definition
\defref{regular integrable systems of (k,n)}, $\Psi$ is transversal
to the submanifold $\sigma(U)$ and is a local diffeomorphism. Hence
$\mathcal{S}=\Psi^{-1}\big(\sigma(U)\big)$ is an embedded submanifold
of $U\times\mathbb{R}^{n}$ with the corank equal to that of $\sigma(U)$,
which is $n$. Note that $\mathcal{S}$ is also a subspace of $\mathcal{L}$
(and therefore $\mathcal{L}$ is a submanifold of $\mathcal{O}_{M}\times\mathbb{R}^{n}$
due to the arbitrariness of $U$). Moreover, $\mathcal{S}$ is closed
in $U\times\mathbb{R}^{n}$. To see this, suppose that $(u_{i},T_{i})\rightarrow(u,T)$
with $(u_{i},T_{i})\in\mathcal{S}$ and $(u,T)\in U\times\mathbb{R}^{n}$.
Then $\Psi(u,T)=\lim\Psi(u_{i},T)=\lim\sigma(u_{i})=\sigma(u)\in\mathcal{S}$.

For any $x\in\mathcal{S}$, the tangent map $\Psi_{*}$ maps $T_{x}\mathcal{S}$
to $T_{\Psi(x)}\sigma(U)$ and therefore $T_{x}\mathcal{S}\cap T_{x}\mathbb{R}^{n}=\{0\}$,
implying $T_{x}\mathcal{S}\oplus T_{x}\mathbb{R}^{n}=T_{x}(U\times\mathbb{R}^{n})$.
As a result, $\pi_{\mathcal{S}}=pr_{U}\big|_{\mathcal{S}}:\mathcal{S}\rightarrow U$
is a local diffeomorphism with $pr_{U}$ : $U\times\mathbb{R}^{n}\rightarrow U$
being the canonical projection. Note that $U\times\{0\}$ is closed
in $U\times\mathbb{R}^{n}$ and therefore is closed in $\mathcal{S}$,
while it is also open in $\mathcal{S}$. Hence $U\times\{0\}$ is
a connected component of $\mathcal{S}$. 

Let $z_{1}$,...,$z_{n}$ be a basis of the the lattice $\mathcal{L}_{o}$.
Shrink $U$ if necessary. Then there exist local sections $\tau_{i}:U\rightarrow\mathcal{S}$
for $i=1,...,n$ such that $\tau_{i}(o)=z_{i}$. Note that $\{z_{1},...,z_{n}\}$
is also a basis of the linear space $\{o\}\times\mathbb{R}^{n}$,
and that for any $w\in\{o\}\times\mathbb{R}^{n}$, $w\in\mathcal{L}_{o}$
if and only if $w=a_{1}z_{1}+\cdots+a_{n}z_{n}$ with $a_{i}$ integers.
Shrinking $U$ again if necessary, one can make the determinant of
the matrix $[\tau_{1}(c),...,\tau_{n}(c)]$ non-zero for each $c\in U$,
and then $\{\tau_{1}(c),...,\tau_{n}(c)\}$ is a basis for $\{c\}\times\mathbb{R}^{n}$. 

Now we show that for each $c\in U$, $\{\tau_{1}(c),...,\tau_{n}(c)\}$
is also a basis of the $\mathbb{Z}^{n}$ lattice $\mathcal{L}_{c}$.
Let $\{z_{1}^{c},...,z_{n}^{c}\}$ be a basis of $\mathcal{L}_{c}$.
Then there exist intergers $k_{1}^{j},...,k_{n}^{j}$ such that $\tau_{j}(c)=k_{1}^{j}z_{1}^{c}+\cdots+k_{n}^{j}z_{n}^{c}$.
Taking the inverse of the matrix $[k_{i}^{j}]$ shows that for each
$w\in\mathcal{L}_{c}$, there exist $p_{1},...,p_{n}\in\mathbb{Q}$
such that $w=p_{1}\tau_{1}(c)+\cdots+p_{n}\tau_{n}(c)$. Hence it
remains to show that for each $w$, the corresponding coeffients $p_{i}=\frac{k_{i}}{m_{i}}$
are integers. 

Denote by $\mathcal{S}_{w}$ the component of $\mathcal{S}$ that
contains $w$. Note that there exits some nonzero integer $m$ such
that $m\cdot p_{i}$ are integers for all $i$. By $*m$ denote the
map $(c,T)\mapsto(c,m\cdot T)$. Then $*m$ is a diffeomorphism on
$\mathcal{O}_{M}\times\mathbb{R}^{n}$ and maps $\mathcal{S}$ to
$\mathcal{S}$. Hence it is a local diffeomorphism on $\mathcal{S}$.
Moreover, as shown above, $\mathcal{S}$ is closed in $U\times\mathbb{R}^{n}$
and therefore $*m$ is a closed map on $\mathcal{S}$. Hence $*m(\mathcal{S}_{w})$
is another connected component of $\mathcal{S}$ which contains $m\cdot w$.
Meanwhile, $\tau:U\rightarrow U\times\mathbb{R}^{n}$: $u\rightarrow m\cdot\big(p_{1}\tau_{1}(u)+\cdots+p_{n}\tau_{n}(u)\big)$
is also a section of $\mathcal{S}\rightarrow U$ and hence $\tau(U)$
is a component of $\mathcal{S}$ which also contains $m\cdot w$.
As a result, $\tau(U)=*m(\mathcal{S}_{w})$. In particular, $m\cdot p_{1}\tau_{1}(o)+\cdots+m\cdot p_{n}\tau_{n}(o)=\tau(o)=m\cdot w_{0}\in\mathcal{L}_{o}$
for some $w_{0}\in\mathcal{S}_{w}\cap\mathcal{L}_{o}$ and hence $w_{0}=p_{1}\tau_{1}(o)+\cdots+p_{n}\tau_{n}(o)=p_{1}z_{1}+\cdots+p_{n}z_{n}$,
implying the $p_{i}$ to be integers. 

Define $\beta:U\times\mathbb{Z}^{n}\rightarrow\pi_{\mathcal{L}}^{-1}(U)$
with $\rho^{-1}(u;a_{1},...,a_{n})=a_{1}\tau_{1}(u)+\cdots+a_{n}\tau_{n}(u)$.
Then $\beta$ is a local trivialization of $\mathcal{L}$ over $U$.
\end{proof}
\begin{cor}
The lattice bundle $\pi_{\mathcal{L}}:\mathcal{L}\rightarrow\mathcal{O}_{M}$
is a locally trivial smooth $\mathbb{Z}^{n}$ bundle with the transition
group $SL(n,\mathbb{Z})$.
\end{cor}

\begin{proof}
Choose a covering $\{U_{\alpha}\}$ of $\mathcal{O}_{M}$ with sections
$\tau_{i}^{\alpha}:U_{\alpha}\rightarrow\mathcal{L}$ forming a basis
of $\mathcal{L}\big|_{U_{\alpha}}$. Reorder $\{\tau_{i}^{\alpha}\big|i=1,...,n\}$
if necessary to fit the orientation of the vector bundle $\mathcal{O}_{M}\times\mathbb{R}^{n}\rightarrow\mathcal{O}_{M}$,
and then for any pair $\alpha,\alpha'$ with $U_{\alpha}\cap U_{\alpha'}$
being non-empty, $\rho_{\alpha'}\circ\rho_{\alpha}^{-1}$ should be
a linear map on $U_{\alpha}\cap U_{\alpha'}\times\mathbb{Z}^{n}$
preserving the orientation and hence $\rho_{\alpha'}\circ\rho_{\alpha}^{-1}\in SL(n,\mathbb{Z})$.
\end{proof}

\subsection{The fibration structure of $M_{\mathcal{C}}$ determined by the monodromy
map}

In the following we give an illustration of how the monodromy map
$\mathcal{M}$ determines the torus fibration $M_{\mathcal{C}}=q_{\Phi}^{-1}(\mathcal{C})\rightarrow\mathcal{C}$.

Recall that $\mathcal{C}\subset\mathcal{O}_{M}$ is a loop with $M_{\mathcal{C}}=q_{\Phi}^{-1}(\mathcal{C})\rightarrow\mathcal{C}$
being a torus fibration. Let $c_{0}$ be a point in $\mathcal{C}$.
The lift of $\mathcal{C}$ to $\mathcal{L}$ gives rise to a group
isomorphism $\mathcal{M}_{\mathcal{C}}$ of $\mathcal{L}_{c_{0}}$,
which is called the monodromy map. Note that it extends canonically
and uniquely to an isomorphism of $\{c_{0}\}\times\mathbb{R}^{n}$,
which we also denote by $\mathcal{M}$. Fix a parametrization $\gamma:[0,1]\rightarrow\mathcal{C}$
with $\gamma(0)=\gamma(1)=c_{0}$. 

According to Corollary \corref{if q_phi has a global section}, $M_{\mathcal{C}}=q_{\Phi}^{-1}(\mathcal{C})\rightarrow\mathcal{C}$
is isomorphic to $\bar{\pi}_{\mathcal{C}}:\big(\mathcal{C}\times\mathbb{R}^{n}\big)\big/\sim_{\mathcal{L}}\rightarrow\mathcal{C}$
with $\bar{\pi}_{\mathcal{C}}([c,T])=c$. Choose an ordered basis
$\bar{w}=(u_{1},...,u_{n})$ of $\mathcal{L}_{c_{0}}$, and let $\tau_{i}:[0,1]\rightarrow\mathcal{L}$
be the lift of $\gamma$ at $u_{i}$ for each $i=1,...,n$ and $\tau_{i}':[0,1]\rightarrow\mathbb{R}^{n}$
be the map such that $\tau_{i}(s)=(\gamma(s),\tau_{i}'(s))$. Then
$\mathcal{M}_{\mathcal{C}}$ is the map mapping $u_{i}=\tau_{i}(0)$
to $\tau_{i}(1)$. Denote by $\mathcal{M}_{\mathcal{C},\bar{w}}$
the matrix representation of $\mathcal{M}_{\mathcal{C}}$ with respect
to $\bar{w}$, i.e. $\mathcal{M}_{\mathcal{C}}(t_{i}u_{i})=(u_{1},...,u_{n})\cdot\mathcal{M}_{\mathcal{C},\bar{w}}\cdot\left[\begin{array}{c}
t_{1}\\
\vdots\\
t_{n}
\end{array}\right]$. In other words, $\mathcal{M}_{\mathcal{C}}(u_{j})=(u_{1},...,u_{n})\cdot\mathcal{M}_{\mathcal{C},\bar{w}}^{j}$
with $\mathcal{M}_{\mathcal{C},\bar{w}}^{j}$ the $j-$th column of
$\mathcal{M}_{\mathcal{C},\bar{w}}$, and hence $(\tau_{1}'(1),...,\tau_{n}'(1))=(\tau_{1}'(0),...,\tau_{n}'(0))\cdot\mathcal{M}_{\mathcal{C},\bar{w}}$. 

Denote by $\sim_{\gamma}$ the equivalence relation on $[0,1]\times\mathbb{R}^{n}$
that identifies $(0,T)$ with $(1,T)$ and by $q_{\gamma}$ the corresponding
quotient map from $[0,1]\times\mathbb{R}^{n}$ to $\mathcal{C}\times\mathbb{R}^{n}$.
Define $r_{\bar{w}}:[0,1]\times\mathbb{R}^{n}\rightarrow[0,1]\times\mathbb{R}^{n}$
by
\begin{equation}
r_{\bar{w}}(s;\left[\begin{array}{c}
t_{1}\\
\vdots\\
t_{n}
\end{array}\right])=\big(s,t_{1}\tau_{1}'(s)+\cdots+t_{n}\tau_{n}'(s)\big).\label{eq:r-def}
\end{equation}
 This is an isomorphism mapping $[0,1]\times\mathbb{Z}^{n}\rightarrow\gamma^{*}(\mathcal{L})$.
Hence the map $h_{\bar{w}}=q_{\gamma}\circ r_{\bar{w}}:[0,1]\times\mathbb{R}^{n}\rightarrow\mathcal{C}\times\mathbb{R}^{n}$
is a bundle epimorphism mapping $[0,1]\times\mathbb{Z}^{n}$ to $\mathcal{L}\big|_{\mathcal{C}}$
and it is a fiberwise isomorphism. Moreover, $h_{\bar{w}}$ identifies
$\{0\}\times\left[\begin{array}{c}
t_{1}\\
\vdots\\
t_{n}
\end{array}\right]$ with $\{1\}\times\mathcal{M}_{\mathcal{C},\bar{w}}^{-1}\cdot\left[\begin{array}{c}
t_{1}\\
\vdots\\
t_{n}
\end{array}\right]$. To see this, first note that since $q_{\gamma}$ identifies $(0,T)$
with $(1,T)$ for every $T\in\mathbb{R}^{n}$, $q_{\gamma}\circ r_{\bar{w}}$
identifies $(0,T)$ with $(1,T')$ if and only if $r_{\bar{w}}(0,T)=(0,T'')$
and $r_{\bar{w}}(1,T)=(1,T'')$. Then with $T=\left[\begin{array}{c}
t_{1}\\
\vdots\\
t_{n}
\end{array}\right]$ and $T'=\left[\begin{array}{c}
t'_{1}\\
\vdots\\
t'_{n}
\end{array}\right]$, according to (\ref{eq:r-def}), it holds that
\[
t_{1}\tau_{1}'(0)+\cdots+t_{n}\tau_{n}'(0)=t'_{1}\tau_{1}'(1)+\cdots+t'_{n}\tau_{n}'(1)=\big(\tau_{1}'(0),...,\tau_{n}'(0)\big)\cdot\mathcal{M}_{\mathcal{C},\bar{w}}\cdot\left[\begin{array}{c}
t'_{1}\\
\vdots\\
t'_{n}
\end{array}\right].
\]
Hence $\left[\begin{array}{c}
t_{1}\\
\vdots\\
t_{n}
\end{array}\right]=\mathcal{M}_{\mathcal{C},\bar{w}}\cdot\left[\begin{array}{c}
t'_{1}\\
\vdots\\
t'_{n}
\end{array}\right]$ and thus $\left[\begin{array}{c}
t'_{1}\\
\vdots\\
t'_{n}
\end{array}\right]=\mathcal{M}_{\mathcal{C},\bar{w}}^{-1}\cdot\left[\begin{array}{c}
t_{1}\\
\vdots\\
t_{n}
\end{array}\right]$.

Note that $quo=q_{\mathcal{L}}\circ h_{\bar{w}}$ is a quotient map
from $[0,1]\times\mathbb{R}^{n}$ to $M_{\mathcal{C}}$ and it charaterizes
the structure of the fibration $M_{\mathcal{C}}\rightarrow\mathcal{C}$.
Denote by $\sim$ the equivalence relation induced by $quo$ on $[0,1]\times\mathbb{R}^{n}$.
$M_{\mathcal{C}}\rightarrow\mathcal{C}$ is isomorphic to $\pi_{\sim}:\big([0,1]\times\mathbb{R}^{n}\big)\big/\sim\rightarrow\mathcal{C}$
with $\pi_{\sim}([s,T])=\gamma(s)$. It can be checked that for each
$s\in[0,1]$, $(s,T)\sim(s,T')$ if and only if $T-T'\in\mathbb{Z}^{n}$.
Then $quo$ (re)factors as $q_{\bar{\eta}}\circ q_{\mathbb{Z}^{n}}$
with $q_{\mathbb{Z}^{n}}:[0,1]\times\mathbb{R}^{n}\rightarrow[0,1]\times\mathbb{T}^{n}$
being the quotient map which sends $(s;t_{1},...,t_{n})$ to $(s;e^{i2\pi t_{1}},...,e^{i2\pi t_{n}})$,
and $q_{\bar{\eta}}$: $[0,1]\times\mathbb{T}^{n}\rightarrow M_{\mathcal{C}}$
being a bundle morphism which is a fiberwise isomorphism. Denote by
$\eta$ the isomorphism $(1,T)\rightarrow(h_{\bar{w}}\big|_{\{0\}\times\mathbb{R}^{n}})^{-1}\circ(h_{\bar{w}}\big|_{\{1\}\times\mathbb{R}^{n}})(1,T)=(0,\mathcal{M}_{\mathcal{C},\bar{w}}\cdot T)$
and by $\bar{\eta}$ the isomorphism $(1,z)\rightarrow(q_{\bar{\eta}}\big|_{\{0\}\times\mathbb{T}^{n}})^{-1}\circ(q_{\bar{\eta}}\big|_{\{1\}\times\mathbb{T}^{n}})(1,z)$.

Then $M_{\mathcal{C}}\cong\big([0,1]\times\mathbb{T}^{n}\big)\big/\sim_{\bar{\eta}}$.
Moreover, the following identity holds
\[
\bar{\eta}\circ q_{\mathbb{Z}^{n}}=q_{\mathbb{Z}^{n}}\circ\eta.
\]

\printbibliography

\end{document}